\documentstyle[11pt]{article}
\input{amssym.def}
\input epsf

\newtheorem{thm}{Theorem}[section]
\newtheorem{lemma}[thm]{Lemma}
\newtheorem{prop}[thm]{Proposition}
\newtheorem{conj}[thm]{Conjecture}
\newtheorem{defn}[thm]{Definition}

\newtheorem{cor}[thm]{Corollary}
\newtheorem{rmk}[thm]{Remark}

\def\index{{\rm{index}}}
\def\codim{{\rm{codim}}}

\def\ker{\mathop{\rm Ker}\nolimits}

\def\im{\mathop{\rm Im}\nolimits \,}
\def\re{\mathop{\rm Re}\nolimits \,}

\def\less{{\prec}}
\def\supp{{\mbox{\rm supp}}}

\newcommand{\C}{{\Bbb C}}

\newcommand{\N}{{\Bbb N}}
\newcommand{\R}{{\Bbb R}}

\newcommand{\cS}{{\cal S}}

\newcommand{\cC}{{\cal C}}

\newcommand{\cP}{{\cal P}}

\newcommand{\cU}{{\cal U}}
\newcommand{\cX}{{\cal X}}

\newcommand{\cV}{{\cal V}}

\begin{document}
\bibliographystyle{plain}

\begin{center}
{\LARGE Gradient-Like Flows and Self-Indexing}

\vspace{.08in}

{\LARGE in Stratified Morse Theory}

\vspace{.23in}

\large{Mikhail Grinberg\footnotemark}

\vspace{.15in}

\large{August 15, 2000}
\end{center}

\footnotetext{Research partially supported by NSF grant
$\#$ DMS-9971030.}

\vspace{.2in}

\begin{abstract}
We develop the idea of self-indexing and the technology of
gradient-like vector fields in the setting of Morse theory on a
complex algebraic stratification.  Our main result is the local
existence, near a Morse critical point, of gradient-like vector
fields satisfying certain ``stratified dimension bounds up to
fuzz'' for the ascending and descending sets.  As a global
consequence of this, we derive the existence of self-indexing
Morse functions.
\end{abstract}

\vspace{.3in}

\tableofcontents

\vfill

\section{Introduction}

The goal of this paper is to develop the idea of self-indexing and
the technology of gradient-like vector fields in the setting of Morse
theory on a complex algebraic stratification.  Our main result
(Theorem \ref{main}) is the local existence, near a Morse critical
point, of gradient-like vector fields satisfying certain ``stratified
dimension bounds up to fuzz'' for the ascending and descending sets.
As a global consequence of this, we derive the existence of
self-indexing Morse functions (Theorem \ref{siinsm}).

This paper traces its roots to the informal lecture notes
``Intersection Homology and Perverse Sheaves'' by R. MacPherson
\cite{MP}.  These notes outline a vision for developing the theory of
middle perversity perverse sheaves on a stratified complex variety,
starting with a definition of a perverse sheaf in the style of the
Eilenberg-Steenrod axioms.  This approach to perverse sheaves relies
on self-indexing Morse functions as a key technical tool.  We now
proceed to introduce this circle of ideas at a leisurely pace.

\subsection{Classical Morse Theory}

Let $X^d$ be a compact smooth manifold, and let $f: X \to \R$ be a
Morse function.  Write $\Sigma_f$ for the set of critical points of
$X$.  The function $f$ is called self-indexing if $f(p) =
\index_f (p)$ for every $p \in \Sigma_f$.  The significance of this
definition is the following.  A self-indexing $f$ gives rise
canonically to a cochain complex $\cC_f = (C^*, d)$ of $\C$-vector
spaces such that:

(a)  $C^i$ has a natural basis up to sign parameterized by
$\{ \, p \in \Sigma_f \; | \; \index_f (p) = i \, \}$;

(b)  $H^*(\cC_f) = H^*(X; \C)$.

Often, one uses an arbitrary Morse function plus another piece
of structure (e.g., a generic metric) to define a complex  for
computing $H^*(X; \C)$.  However, one needs self-indexing to produce
a canonical complex starting with a Morse function alone.  The
complex $\cC_f$ is called the Morse-Smale complex.  It plays a
central role in Milnor's beautiful exposition \cite{Mi} of Smale's
proof of the $h$-cobordism theorem in dimensions $> 4$.  The first
step of that proof is the following theorem.

\begin{thm}\label{siincm} {\em [16, Theorem B]}
Every compact, smooth  manifold $X$ admits a self-indexing Morse
function $f: X \to \R$.
\end{thm}

We recall an outline of the proof of this in \cite{Mi}.  Start with
any Morse function $g : X \to \R$.  A $\nabla g$-like vector field
$V$ on $X$ is a $C^\infty$ vector field such that:

(a)  $V_p = 0$ for all $p \in \Sigma_g$;

(b)  $V_x \, g > 0$ for all $x \notin \Sigma_g$.

\noindent
(This definition is slightly different from the one in \cite{Mi}.)
Given such a $V$, let $\psi_V : X \times \R \to X$ denote its flow.
Then for $p \in \Sigma_g$, let
$$M^\pm_V (p) = \{ \, x \in X \; | \; \lim_{t \to \mp \infty}
\psi_V (x,t) = p \, \}.$$
The sets $M^\pm_V (p)$ are called the ascending and descending sets
of $p$ (relative to $V$).  Here now are the three main steps of the
proof of Theorem \ref{siincm}.

{\em Step 1: Morse lemma.}  By Morse lemma, there exists a
$\nabla g$-like vector field $V$ such that, for every $p \in \Sigma_g$,
the set $M^-_V (p)$ is a manifold of dimension $\index_g(p)$, and
$M^+_V (p)$ is a manifold of dimension $d - \index_g(p)$.

{\em Step 2: General position.}  By perturbing the vector field $V$,
we can ensure that $M^-_V (p) \cap M^+_V (q) = \emptyset$ whenever
$\index_g (p) \leq \index_g(q)$. 

{\em Step 3: Modifying $g$.}  There exists a Morse function
$f : X \to \R$ such that:

(a)  near each $p \in \Sigma_g$, we have $f(x) = g(x) + \index_g(p)
- g(p)$;

(b)  $V$ is $\nabla f$-like.

\noindent
Conditions (a) and (b) guarantee that $f$ is self-indexing.  The role
of the vector field $V$ is to ensure that $\Sigma_f = \Sigma_g$.

The main result of this paper (Theorem \ref{main}) is a replacement
for Step 1 in the complex stratified setting.  Steps 2 and 3 carry
over in a more or less straightforward fashion.

\subsection{Stratified Morse Theory}

Let now $X^d$ be a smooth complex algebraic variety, and let $\cS$ be
an algebraic Whitney stratification of $X$.  Stratified Morse theory,
pioneered by Goresky and MacPherson (see \cite{GM1}, \cite{GM2}), aims
to study the topology of the pair $(X, \cS)$ by means of a real
$C^\infty$ function $f: X \to \R$.  The following definitions go back
to the original paper \cite{GM1}.

\begin{defn}
(i)   Let $p \in X$ be a point contained in a stratum $S$.  We say
that $p$ is critical for $f$ ($p \in \Sigma_f$) if it is critical for
the restriction $f|^{}_S$.

(ii)  For $S \in \cS$, let $\Lambda_S$ be the conormal bundle
$T^*_S X \subset T^*X$, and let $\Lambda = \Lambda_\cS = \bigcup_S
\Lambda_S$.  This $\Lambda$ is called the conormal variety to $\cS$;
it is a closed Lagrangian subvariety of $T^*X$ (the fact that
$\Lambda$ is closed is equivalent to one of the Whitney conditions).
Note that $p \in \Sigma_f$ if and only if $d_p f \in \Lambda$. 

(iii) Let $p \in \Sigma_f$.  We say that $p$ is Morse for $f$ if it
is Morse for the restriction $f|^{}_S$ (where $S$ is the stratum
containing $p$) and $d_p f \in \Lambda^0$, the smooth part of $\Lambda$.
The set $\Lambda^0$ is called the set of generic conormal vectors to
$\cS$.
\end{defn}

It is important to note that there is nothing specifically complex
about the above definitions: they would apply equally well to any
Whitney stratification of a real $C^\infty$ manifold.  By contrast,
the following definition is only justified in the complex algebraic
(or analytic) setting.

\begin{defn}
Let $p \in \Sigma_f$ be Morse and let $S$ be the stratum of $p$.
We define
$$\index_f (p) = \index_{f|^{}_S} (p) - \dim^{}_\C S = 
\index_{f|^{}_S} (p) + \codim^{}_\C S - d.$$
\end{defn}

The normalization constant $d = \dim_\C X$ is subtracted to make the
possible range of values of the index symmetric about the origin: 
$\index_f (p) \in \{ -d, \dots, d \}$.  Ignoring this normalization,
the index is attempting to count the (real) ``descending directions''
of $f$ at $p$, by counting the tangent directions in the obvious way
and then postulating that exactly half of the normal directions are
descending.  The point of this paper is that this method of counting
is justified.  The justification comes from a local result concerning
Morse theory near a point stratum.  Before stating the result itself,
we present a somewhat sharper conjecture.

\begin{conj}\label{mainconj}
Let $X = \C^d$, $p \in X$ be the origin, and $\cS$ be an algebraic
Whitney stratification of $X$ such that $\{p\}$ is a stratum.  Let
$\Lambda^0_p = \Lambda_{ \{ p \} } \cap \Lambda^0$ be the set of
generic covectors at $p$, and let $f \in \Lambda^0_p$.  Regard $f$
as a linear function $f : X \to \R$.  Then there exists a closed
ball $B$ around $p$ and an $\cS$-preserving $\nabla f$-like vector
field $V$ on $B$ such that the ascending and descending sets
$M^\pm_V(p)$ satisfy:
$$\dim^{}_\R M^\pm_V(p) \cap S \leq \dim^{}_\C S \;\;\;
\mbox{for every} \;\;\; S \in \cS.$$
\end{conj}

The reader is referred to Section 2.2 for precise definitions of an
$\cS$-preserving $\nabla f$-like vector field and of the sets
$M^\pm_V(p)$.  However, the truth and the level of difficulty of this
conjecture might be sensitive to the exact class of vector fields
chosen.  Therefore, it is best to view the conjecture as being
somewhat imprecise, with the phrase ``$\cS$-preserving $\nabla f$-like
vector field'' open to interpretation.  Theorem \ref{main} should be
seen as Conjecture \ref{mainconj} ``up to fuzz''.  

\begin{thm}\label{main}
Let $X$, $p$, $\cS$ be as in Conjecture \ref{mainconj}, and let
$\Delta \subset \Lambda^0_p$ be an open set.  Then there exit an
$f \in \Delta$ (which we view as a linear function $f : X \to \R$),
a closed ball $B$ around $p$, and a closed, real semi-algebraic
$K \subset B$ such that:

(i)   $\dim^{}_\R K \cap S \leq \dim^{}_\C S$ for every $S \in \cS$;

(ii)  $f^{-1}(0) \cap K = \{ p \}$; and

(iii) for every open $\cU \supset K$, there exists an $\cS$-preserving
$\nabla f$-like vector field $V$ on $B$ with $M^\pm_V(p) \subset \cU$.
\end{thm}

Theorem \ref{main} will be proved in Sections 4 and 5.
It is important to note that for both Conjecture \ref{mainconj}
and Theorem \ref{main} it is crucial that the pair $(X, \cS)$ is
assumed to be complex.  For contrast, consider $X = \R^3$, let $\cS$
be the obvious stratification with six strata whose $2$-skeleton is
the cone $x^2+y^2=z^2$, and take $f(x,y,z) = z$.  Then every
$\cS$-preserving $\nabla f$-like vector field $V$ on $X$ will have
$\dim M^+_V(p) = \dim M^-_V(p) = \dim X = 3$.  (An ice cream cone can
hold a 3d amount of ice cream!)  This example shows that the concept
of index does not generalize well to Morse theory on a real
stratification.  It is also easy to modify the above example to one
with all strata of even dimension (just put the same cone in $\R^4$).
The following consequence of Theorem \ref{main} will be proved in
Section 6.

\begin{thm}\label{siinsm}
Every proper, Whitney stratified complex algebraic variety $(X, \cS)$
admits a self-indexing Morse function $f : X \to \R$.
\end{thm}

\subsection{Perverse Sheaves}

We now explain the motivation behind Theorems \ref{main} and
\ref{siinsm} coming from the theory of middle perversity perverse
sheaves.  Let $(X^d, \cS)$ be as in Section 1.2.  Associated to the
pair $(X, \cS)$ is the category $\cP(X, \cS)$ of middle perversity
perverse sheaves on $X$ constructible with respect to $\cS$.  This
category was first introduced by Beilinson, Bernstein, and Deligne in
the seminal paper \cite{BBD} in 1982.  They defined $\cP(X, \cS)$ as a
full subcategory of $D^b_\cS (X)$, the bounded $\cS$-constructible
derived category of sheaves on $X$.  This way $\cP(X, \cS)$ inherits
an additive structure from $D^b_\cS (X)$.  In fact, it turns out to
be an abelian category, unlike the derived category $D^b_\cS (X)$.
Perverse sheaves are ``simpler'' than the derived category in several
other respects.  For example, they form a stack (i.e., objects and
morphisms are locally defined), which is not the case for
$D^b_\cS (X)$.

For these reasons, it was long felt desirable to have a definition
of $\cP(X, \cS)$ which does not rely on the more complicated
derived category, and which elucidates the very simple formal
properties of perverse sheaves.  Such a definition was proposed by
MacPherson in his AMS lectures in San Francisco in January 1991.
Informal notes of those lectures \cite{MP} produced by MacPherson (and
dated 1990) were distributed at the meeting, but were never published.
A slightly modified version of the same definition has now appeared
in print in the lecture notes \cite{GrM} from the summer institute
held in Park City, Utah in July 1997.

Roughly, MacPherson's definition (the 1997 version) proceeds as
follows.  A perverse sheaf $P \in \cP(X, \cS)$ is an assignment
$(Y,Z) \mapsto H^*(Y, Z; P)$, plus coboundary and pull-back maps,
subject to a set of axioms modeled after the Eilenberg-Steenrod
axioms.  Here $(Y,Z)$ is a {\em standard pair} in $X$, that is a
pair of compact subsets of a certain special kind.  Namely, $Y$ is
a smooth, real $2d$ dimensional submanifold of $X$ with corners of
codimension two, and $Z$ is a union of some of the walls of $Y$,
so that each corner is formed by one wall from $Z$ and one wall
not from $Z$, and all the boundary strata of $Y$ are transverse
to $\cS$.  It is a very pleasant feature of this definition that we
only need to consider standard pairs.  The structure group
$H^*(Y, Z; P)$ is called the cohomology of $Y$ relative to $Z$ with
coefficients in $P$; it is a finite dimensional, graded vector space
over $\C$.

Turning now to the axioms, all but one (the dimension axiom)
represent a more or less straightforward adaptation of the
classical Eilenberg-Steenrod to the setting of pairs inside a
fixed space.  The remaining axiom runs as follows.

\vspace{.1in}

\noindent
{\bf Dimension Axiom:}  Let $f : X \to \R$ be a proper Morse
function, and $c \in \R$ be a value such that $f^{-1}(c)$
contains exactly one critical point: $f^{-1}(c) \cap \Sigma_f =
\{ p \}$.  Then
$$H^i (\{ x \in X \; | \; c-\epsilon \leq f(x) \leq c+\epsilon\}, \,
\{ x \in X \; | \; f(x) = c-\epsilon\}; P) = 0$$
for $i \neq \index_f(p)$ and $0 < \epsilon \ll 1$.

\vspace{1in}

Here is how this definition fits with the idea of self-indexing.
Given a standard pair $(Y,Z)$, one can use a variant of Theorem
\ref{siinsm} to find a self-indexing function $$f: Y \setminus
\{ \mbox{corners} \} \to [-d-1, d+1]$$ which is {\em adapted to}
the pair $(Y,Z)$.  More precisely, this means that
$f^{-1}(-d-1) = Z \setminus \{ \mbox{corners} \}$, 
$f^{-1}(d+1) = \partial Y \setminus Z$, and the level sets of $f$
``look like pages of a book'' near the corners.  Now, just as in
classical Morse theory, such a function  $f$ gives rise to a functor
$$\cC_f : \cP(X, \cS) \to \{ \mbox{cochain complexes} \}$$
with the property that $H^*(\cC_f(P)) = H^*(Y, Z; P)$ for every
$P \in \cP(X, \cS)$.  If we write $\cC_f = (C^*, d)$, then the
functor $$C^i : \cP(X, \cS) \to \{ \mbox{vector spaces} \}$$ breaks
up into a direct sum over $\{ \, p \in \Sigma_f \; | \; \index_f (p)
= i \, \}$ of the so-called {\em Morse group} functors
$$M_{d_p f} : \cP(X, \cS) \to \{ \mbox{vector spaces} \}.$$
Each of the $M_{d_p f}$ is localized near $p$, i.e., factors through
the restriction to any open neighborhood of $p$.  (Perverse sheaves
on open subsets of $X$ and restriction functors between them are
defined in the obvious way.)  Moreover, as implied by the notation,
the functor $M_{d_p f}$ depends on $f$ only through the differential
$d_p f \in \Lambda^0$.

The existence of the lift $\cC_f$ from cohomology to cochain
complexes with the above properties, by itself, has important
consequences for the structure of the category $\cP(X, \cS)$.
For example, it immediately implies that morphisms of perverse
sheaves are locally defined.  More is true however: the functor
$\cC_f$ is, in some sense, independent of $f$; the only essential
dependence is on the pair $(Y,Z)$.  To be precise, given any two 
self-indexing functions $f_1, f_2$ adapted to $(Y,Z)$, the functors
$\cC_{f_1}$ and $\cC_{f_2}$ are related by a quasi-isomorphism which
is itself canonical up to chain homotopy.  A reader familiar with
Floer homology has, no doubt, anticipated this.  However, unlike in
Floer homology, this independence of the function requires a proof
which is rather more involved than the proof of $d^2 = 0$.
(The later is a straightforward consequence of the axioms.)

Once the independence of the function is established in a suitably
flexible form, all of the (soft) formal properties of perverse
sheaves can be derived from considering the complexes $\cC_f$ for
different functions and standard pairs.  For example, the abelian
property of $\cP(X, \cS)$ follows easily from the abelian property
of the category of cochain complexes.  Finally, one can develop the
theory of perverse sheaves, starting with MacPherson's definition,
to the point of proving the following.

\begin{thm}\label{indep}
MacPherson's definition of $\cP(X, \cS)$ agrees with the original
definition due to Beilinson-Bernstein-Deligne.
\end{thm}

The present author has recently had an opportunity to work through
a proof of Theorem \ref{indep} in a graduate course at MIT, and a
future paper describing this proof is being planned.  The goal of
the present paper is to provide the main geometric tool:
self-indexing Morse functions.  We conclude this introduction with a
conjecture due to MacPherson.

\begin{conj}\label{fary}
Omitting the dimension axiom from MacPherson's definition of
$\cP(X, \cS)$ gives a category which is equivalent to $D^b_\cS (X)$.
\end{conj}

\noindent
{\bf Acknowledgments.}  I would like to thank the IHES in
Bures-sur-Yvette for their hospitality during August 99 when most of
this work was done.  My intellectual debt to Bob MacPherson is clear
from the preceding introduction.  I am also grateful to the students
in my course at MIT for patiently (but not uncritically) allowing me
to think through these ideas out loud.

\section{Technical Preliminaries}

In this section, we summarize the preliminary material on
controlled and weakly controlled vector fields that we will
need, and define precisely the ingredients of Theorem \ref{main}.
The reader is referred to [2, Chapter 2] and [14, Chapter 2.5]
for a detailed treatment of controlled vector fields.  However,
our definitions differ slightly from these sources.  The main
distinction is that our notion of a quasi-distance function is
more flexible than the corresponding notions in \cite{G-al} and
\cite{dPW}.  This is necessary to make Theorem \ref{isotopy} true
``in the controlled category.''

\subsection{Controlled Vector Fields}

Let $X$ be a real $C^\infty$ manifold with a Whitney
stratification $\cS$.

\begin{defn}\label{proj}
Let $S \in \cS$ be a stratum, and let $U_S$ be an open
neighborhood of $S$.  A tubular projection $\Pi_S : U_S \to S$
is a smooth submersion restricting to the identity on $S$.
\end{defn}

\begin{defn}\label{qdprelim}
Let $M$ be a smooth manifold, and let $E \to M$ be a vector
bundle with zero section $Z$.  A quasi-norm on $E$ is a smooth
function $\rho : E \setminus Z \to \R_+$ such that
$\rho(\lambda \, e) = \lambda \, \rho(e)$ for every
$e \in E \setminus Z$ and $\lambda \in \R_+$.
\end{defn}

\begin{defn}\label{qd}
Let $S \in \cS$ be a stratum, and let $U_S$ be an open
neighborhood of $S$.  \\ A quasi-distance function
$\rho : U_S \setminus S \to \R_+$ is a smooth function
satisfying the following condition.  There exist a vector
bundle $p : E \to S$ with zero section $Z$, an open neighborhood
$U' \subset E$ of $Z$, and a diffeomorphism $\phi : U' \to U$,
such that $\phi|^{}_Z = p|^{}_Z$ and $\rho \circ \phi$ is the
restriction to $U' \setminus Z$ of a quasi-norm on $E$.
\end{defn}

\begin{defn}\label{cdat}
Control data on $(X, \cS)$ is a collection $\{ U_S, \Pi_S, \rho_S
\}_{S \in \cS}$, where $U_S \supset S$ is an open neighborhood,
$\Pi_S : U_S \to S$ is a tubular projection, and $\rho_S : U_S
\setminus S \to \R_+$ is a quasi-distance function, subject to the
compatibility conditions $\Pi_S \circ \Pi_T = \Pi_S$ and $\rho_S
\circ \Pi_T = \rho_S$, both being equalities between maps of $U_S
\cap U_T$.
\end{defn}

\begin{defn}\label{compat}
Let $\{ U_S, \Pi_S, \rho_S \}$ be control data on $(X, \cS)$,
let $\cU \subset X$ be an open subset, let $A$ be a set,
and let $f : \cU \to A$ be a map of sets.  We say that
$\{ U_S, \Pi_S, \rho_S \}$ is $f$-compatible on $\cU$ if, for
every $S \in \cS$, there is a neighborhood $U'_S$ of $S \cap \cU$
such that $f \circ \Pi_S = f$ on $U'$.
\end{defn}

\begin{defn}\label{chom}
Let $(X, \cS)$, $(\hat X, \hat \cS)$ be two Whitney stratified
manifolds, with control data $\{ U_S, \Pi_S, \rho_S \}$ and
$\{ U_{\hat S}, \Pi_{\hat S}, \rho_{\hat S} \}$.  A controlled
homeomorphism $\phi : X \to \hat X$ is a homeomorphism which takes
strata diffeomorphicly onto strata, establishing a bijection
$S \mapsto \hat S$, and satisfies the following condition.  For
every $S \in S$, there is a neighborhood $U'_S \subset U_S$ of $S$
such that $\phi \circ \Pi_S = \Pi_{\hat S} \circ \phi$ and
$\phi \circ \rho_S = \rho_{\hat S} \circ \phi$ on $U'_S$.
\end{defn}

\begin{defn}\label{cvf}
Suppose we are given control data $\{ U_S, \Pi_S, \rho_S\}$ on
$(X, \cS)$.  A controlled vector field $V$ on $X$ compatible with
$\{ U_S, \Pi_S, \rho_S\}$ is a collection $\{ V_S \}_{S \in \cS}$
of smooth vector fields on the individual strata, satisfying the
following condition.  For every $S \in \cS$, there exists a
neighborhood $U'_S \subset U_S$ of $S$ such that:

(a) $(\Pi_S)_* V_x = V_{\Pi_S (x)}$ for every $x \in U'_S$;

(b) $V_x \, \rho_S = 0$ for every $x \in U'_S \setminus S$.
\end{defn}

Integrating controlled vector fields is a basic technique for
constructing controlled homeomorphisms, going back to the work
of Thom \cite{Th} and Mather \cite{Ma}.  The following variant of
[15, Lemma 4.11] and [13, Theorem 1.1] will serve as our basic
tool for constructing controlled vector fields.

\begin{lemma}\label{ecvf}
Let $\{ U_S, \Pi_S, \rho_S\}$ be control data on $(X, \cS)$.  Let
$S$ be a stratum, let $U \subset S$ be open in $S$, and let $V$ be
a smooth vector field on $U$.  Then there exist an open $\cU \subset
X$ with $\cU \cap S = U$ and a controlled vector field $\tilde V$ on
$\cU$, compatible with $\{ U_S, \Pi_S, \rho_S\}$, such that
$\tilde V|^{}_U = V$.  Furthermore, the vector field $\tilde V$ can
be chosen to be continuous as a section of $T\cU$. \hfill$\Box$
\end{lemma}

\subsection{Weakly Controlled Vector Fields}

Controlled vector fields are not suitable for discussing the
ascending and descending sets.  Indeed, the trajectory of a
controlled vector field can not approach a point on a smaller
stratum as time tends to infinity.  We will therefore need the
following definition.

\begin{defn}\label{tvf}
Suppose we are given control data $\{ U_S, \Pi_S, \rho_S\}$ on
$(X, \cS)$.  A weakly controlled vector field $V$ on $X$ is a
collection $\{ V_S \}_{S \in \cS}$ of smooth vector fields on the
individual strata, satisfying the following condition.  For every
$S \in \cS$, there exists a neighborhood $U'_S \subset U_S$ of $S$
and a number $k > 0$ such that:

(a) $(\Pi_S)_* V_x = V_{\Pi_S (x)}$ for every $x \in U'_S$;

(b) $|V_x \, \rho_S| < k \cdot \rho_S(x)$ for every
$x \in U'_S \setminus S$.
\end{defn}

\begin{prop}\label{tvf2thom}
Weakly controlled vector fields integrate to stratum preserving
homeomorphisms.  More precisely, let $V$ be a weakly controlled
vector field on $X$.  Then for every $x \in X$, there is a
neighborhood $U_x$ of $x$ and number $t_0 > 0$ such that, for
every $t \in [-t_0, t_0]$, the time-$t$ flow of $V$ is defined
on $U_x$ and gives a stratum preserving homeomorphism
$\psi_{V,t} : U_x \to \psi_{V,t} (U_x)$, which is smooth on each
stratum.
\end{prop}

\noindent
{\bf Proof:}  The main thing to check is that a trajectory of
$V$ lying in a stratum $T$ can not approach a point on a smaller
stratum $S \subset \bar T$ in finite time.  This follows from
condition (b) in Definition \ref{tvf}.  See [14, Proposition 2.5.1]
for more details.
\hfill$\Box$

\vspace{.1in}

The next two definitions clarify the statement of Theorem \ref{main}.

\begin{defn}\label{gradlike}
Let $f : X \to \R$ be a $C^\infty$ function.  An $\cS$-preserving
$\nabla f$-like vector field $V$ on an open subset $\cU \subset X$ 
is a weakly controlled vector field (compatible with some control
data on $(X, \cS)$) satisfying:

(a)  $V_p = 0$ for all $p \in \Sigma_f \cap \cU$;

(b)  $V_x \, g > 0$ for all $x \in \cU \setminus \Sigma_f$.

\noindent
A $\nabla f$-like vector field on an arbitrary subset
$A \subset X$ (e.g., on a closed ball) is the restriction to $A$
of a $\nabla f$-like vector field on some open $\cU \supset A$.
\end{defn}

\begin{defn}\label{andmfds}
Let $f : X \to \R$ be a $C^\infty$ function, let $V$ be a
$\nabla f$-like vector field on some $A \subset X$, and let
$p \in \Sigma_f \cap A$.  We define $M^-_V (p)$ to be the set
of all $x \in A$ such that the trajectory $\psi_{V,t} (x)$ is
contained in $A$ for all $t \geq 0$ and we have
$\displaystyle \lim_{t \to \infty} \psi_{V,t} (x) = p$.  The
ascending set $M^+_V (p)$ is defined similarly.
\end{defn}

\subsection{The Flow Topology}

In this section, we discuss the notion of the flow topology on
the set of weakly controlled vector fields.  It will give us a
degree of flexibility, making some of our constructions less
tied to the choice of control data.

Let $\cV(X, \cS)$ be the set of all weakly controlled vector
fields on $(X,\cS)$, compatible with all possible control data,
and let $\cV(X)$ be the union of the $\cV(X, \cS)$ over all
Whitney stratifications $\cS$ of $X$.  Fix a Riemannian metric
$g$ on $X$.  Let $V \in \cV(X)$, let $K \subset X$ be a compact
subset, let $t > 0$ be a number such that $\psi_{V,s}(x)$
is defined for all $x \in K$ and all $s \in [-t, t]$, and let
$\epsilon > 0$ be any positive number.  Define
$$\cU(V, K, t, \epsilon) = \{ V' \in \cV \; | \;
\forall x \in K, s \in [-t, t] : \mbox{dist}_g
(\psi_{V,s}(x), \, \psi_{V',s}(x)) < \epsilon \},$$
where we set $\mbox{dist}_g (\psi_{V,s}(x), \, \psi_{V',s}(x)) = 
+ \infty$ if $\psi_{V',s}(x)$ is undefined.

\begin{defn}\label{flow-t}
The flow topology on $\cV(X)$ is the weakest topology in which
all the $\cU(V, K, t, \epsilon)$ are open sets.
\end{defn}

It is easy to see that the flow topology is independent of the
metric $g$.

\begin{prop}\label{approx}
(i)  Let $\cS$ and $\hat \cS$ be two Whitney stratifications of $X$
such that $\hat \cS$ is a refinement of $\cS$.  Fix a vector field
$\hat V \in \cV(X, \hat \cS)$ and control data $\{ U_S, \Pi_S,
\rho_S \}$ on $(X, \cS)$.  Then there exists a sequence $\{ V_i
\in \cV(X, \cS) \}_{i \in \N}$ of vector fields compatible with
$\{ U_S, \Pi_S, \rho_S \}$ such that $V_i \to \hat V$ in the flow
topology.

(ii) In the situation of part (i), assume that $\hat \cS$ has a unique
point stratum $\{ p \}$, and that $f : X \to \R$ is a smooth function
whose only stratified critical point with respect to $\hat \cS$ is $p$.
Assume also that $V$ is $\nabla f$-like.  Then all the $V_i$ can be
chosen to be $\nabla f$-like too.
\end{prop}

\noindent
{\bf Proof:}  This is an exercise using Lemma \ref{ecvf} and
partitions of unity.  For part (i), the first step is to show that
there is a sequence $\{ \hat V_i \in \cV(X, \hat \cS)\}$ of continuous
vector fields compatible with the same control data as $V$, such that
$\hat V_i \to \hat V$ in the flow topology.  The second step is to
show that each $\hat V_i$ can be approximated in the $C^0$ topology
by continuous vector fields from $\cV(X, \cS)$ compatible with
$\{ U_S, \Pi_S, \rho_S \}$.  Then it remains to note that the $C^0$
topology on continuous weakly controlled vector fields is stronger
than the flow topology.  Part (ii) is similar.
\hfill$\Box$

\section{Stratified Morse Lemma}

In this section, we prove an isotopy lemma (Theorem
\ref{isotopy}) adapted to the local study of stratified Morse 
functions.  As consequences, we derive some local normal form
statements, one of which (Corollary \ref{ml}) may be called
the stratified Morse lemma.  The results of this section are
very close to those of H. King in \cite{K1} and \cite{K2}.  We
continue with a Whitney stratified smooth manifold $(X, \cS)$.
As in Section 1.2, we denote by $\Lambda^0$ the set of generic
conormal vectors to $\cS$.  We also let $\Lambda^0_S = \Lambda^0
\cap \Lambda_S$, for each $S \in \cS$.

\begin{thm}\label{isotopy}
Let $S \in \cS$ with $\dim S = s$.  Let $B \subset S$ be
a closed $s$-ball smoothly embedded in $S$, with interior
$B^\circ \subset B$ and a fixed point $a \in B^\circ$.  Let
$\cU \subset X$ be an open neighborhood of $B$, let $p : \cU
\to S$ be a smooth submersion restricting to the identity on
$\cU \cap S$, and let $f: \cU \to \R$ be a smooth function
such that $f|^{}_{S \cap \cU} = 0$ and $d_b f \in \Lambda^0_S$ 
for every $b \in B$.  Write $N = p^{-1} (a)$.  Then there
exist an open set $U \subset \cU$ with $U \cap S = B^\circ$,
control data on $U$ which is $p$-compatible on $U$ and
$f$-compatible on $U \setminus S$, and a controlled
homeomorphism $\phi : (U \cap N) \times B^\circ \to U$
such that:

(i)   $\phi \, (x,a) = x$ for every $x \in U \cap N$;

(ii)  $p \circ \phi \, (x,y) = y$ for every
$x \in U \cap N$ and $y \in B^\circ$.

(iii) $f \circ \phi \, (x,y) = f(x)$ for every
$x \in U \cap N$ and $y \in B^\circ$.
\end{thm}

\begin{cor}\label{ml}
Let $f : X \to \R$ be a smooth function with a stratified
Morse critical point $a$, lying in a stratum $S$, let
$\cU \subset X$ be a neighborhood of $a$, and let $p : \cU
\to S$ be a smooth submersion restricting to the identity on
$\cU \cap S$.  Write $N = p^{-1} (a)$.  Then there exist a
smaller neighborhood $U \subset \cU$ of $a$, control data on
$U$ which is $p$-compatible on $U$ and $f$-compatible on
$U \setminus S$, and a controlled homeomorphism
$\phi : (U \cap N) \times (U \cap S) \to U$ such that:

(i)  $\phi \, (x,a) = x$ for every $x \in U \cap N$;

(ii) $p \circ \phi \, (x,y) = y$ for every
$x \in U \cap N$ and $y \in U \cap S$;

(iii) $f \circ \phi \, (x,y) = f(x) + f(y) - f(a)$
for every $x \in U \cap N$ and $y \in U \cap S$. 
\end{cor}

\noindent
{\bf Proof:}  Apply Theorem \ref{isotopy} to the function
$f - f \circ p$.
\hfill$\Box$

\begin{rmk}{\em
Corollary \ref{ml} can be used to give a short proof of Goresky and
MacPherson's product theorem for Morse data, as stated in Chapter 1 of
\cite{GM2} (Theorem STM, part B).  This was earlier observed by H. King
in \cite{K1} and \cite{K2}.  We should also mention that a different
short proof of Goresky and MacPherson's result has appeared recently in
\cite{Ha}.}
\end{rmk}

\begin{cor}\label{ns}
Let $S$ be a stratum.  Suppose we have two points $a, \hat a \in S$,
two normal slices $N, \hat N$ to $S$ passing through these points,
and two functions $f : N \to \R$ and $\hat f : \hat N \to \R$ with
$f(a) = \hat f (\hat a) = 0$.  Assume that the differentials $d_a f$
and $d_{\hat a} \hat f$ are both in $\Lambda^0_S$ and, moreover, in
the same path-component of $\Lambda^0_S$.  Then there exist open
neighborhoods $U \subset N$ and $\hat U \subset \hat N$ of $a$ and
$\hat a$, control data on $U$ and $\hat U$, $f$- (resp. $\hat f$-)
compatible on $U \setminus \{ a \}$ (resp. $\hat U \setminus
\{ \hat a \}$), and a controlled homeomorphism $\phi : U \to \hat U$
such that $f = \hat f \circ \phi$.
\end{cor}

\noindent
{\bf Proof:}  Apply Theorem \ref{isotopy} to a suitable
function on the product $X \times \R$.
\hfill$\Box$

\begin{rmk}{\em
As another corollary of Theorem \ref{isotopy}, we note the fact
that stratified Morse functions are topologically locally stable
in the same sense as the ordinary Morse functions.  In other words,
a stratified Morse function is (locally) right-equivalent, by a
controlled homeomorphism, to any nearby function with the same
critical value.  We omit the precise statement.}
\end{rmk}

The proof of Theorem \ref{isotopy} is based on the following
lemma.

\begin{lemma}\label{clever}  In the situation of Theorem
\ref{isotopy}, fix a Riemannian metric $g$ on $X$.  Let
$r : \cU \to \R_{\geq 0}$ be the distance-to-$S$ function.
Let $\alpha = (f,r): \cU \setminus S \to \R \times \R_+$.
Then there exist an open neighborhood $\cU' \subset \cU$ of $B$ and
number $k > 0$, such that $(p, \alpha) : \cU \setminus S \to
S \times \R \times \R_+$ is a stratified submersion on the set
$\{ x \in \cU' \setminus S \; | \; |f(x)| \leq k \cdot r(x) \}$.
\end{lemma}

\noindent
{\bf Proof:}  Suppose the lemma is false.  Let $\Sigma \subset
\cU \setminus S$ be the stratified critical locus of the map
$(p, \alpha)$.  Then there exists a sequence $\{ x_i \in \Sigma \}$,
converging to a point $b \in B$, such that
\begin{equation}\label{inclever}
\lim_{i \to \infty} \frac{f(x_i)}{r(x_i)} = 0.
\end{equation}
By passing if necessary to a subsequence, we can assume that all
the $\{ x_i \}$ lie in the same stratum $R$, and that there
exists a limit $\Delta = \lim T_{x_i} R \subset T_b X$.  But then,
combining equation (\ref{inclever}) with the Whitney conditions
for the pair $(S, R)$, we may conclude that the differential
$d_b f$ annihilates $\Delta$.  This, however, contradicts the
genericity assumption $d_b f \in \Lambda^0_S$.
\hfill$\Box$

\vspace{.1in}

\noindent
{\bf Proof of Theorem \ref{isotopy}:}  Without loss of generality, we
may assume that $X$ is the total space of a vector bundle $X \to M$,
that $S$ is the zero section, and that $f : \cU \to \R$ is the
restriction of a fiber-wise linear function $\tilde f : X \to \R$.
Fix a Euclidean structure in the bundle $X \to M$, and let
$r : X \to \R_{\geq 0}$ be the associated norm.
Let $\alpha = (f,r): \cU \setminus S \to \R \times \R_+$, and apply
Lemma \ref{clever} to produce a neighborhood $\cU' \subset \cU$ of
$B$ and a number $k > 0$.

{\it Step 1.}  Let $R_k = \{ (\xi, \eta) \in \R \times \R_+ \; |
\; |\xi| \leq k \cdot \eta \}$.  There exists a smooth,
$R_+$-equivariant function $\tilde\rho : \R \times \R_+ \to \R_+$
such that $\tilde\rho(\xi, \eta) = |\xi|$ for all $(\xi, \eta)
\notin R_k$.

{\it Step 2.}  Let $\rho : \cU' \setminus S \to \R_+$ be the
composition $\rho = \tilde\rho \circ \alpha$.  Then $\rho$ is a
quasi-distance function for $S$.

{\it Step 3.}  There exist control data $\{ U_T, \Pi_T, \rho_T \}$
on $\cU'$ such that $U_S = \cU'$, $\Pi_S = p|^{}_{\cU'}$,
$\rho_S = \rho$, and $\{ U_T, \Pi_T, \rho_T \}$ is
$\alpha$-compatible on some open set containing
$\cU' \cap \alpha^{-1} (R_k)$.

{\it Step 4.}  Every smooth vector field $V$ on $B$ extends to a
controlled vector field $\tilde V$ on $\cU' \cap p^{-1} (B)$,
compatible with $\{ U_T, \Pi_T, \rho_T \}$ and satisfying
$d \alpha (\tilde V_x) = 0$ for every $x \in \cU' \cap p^{-1} (B)
\cap \alpha^{-1} (R_k)$.

{\it Step 5.}  The homeomorphism $\phi$ is constructed by
integrating the $\tilde V_i$ for a suitable collection
$\{ V_i \}_{i=1}^s$ of vector fields on $B$.
\hfill$\Box$

\section{Construction of the Set $K$}

We now begin the proof of Theorem \ref{main}.  It is based
on the following lemma.

\begin{lemma}\label{flag}
In the situation of Theorem \ref{main}, there exist a covector
$l \in \Delta$ (which we regard as a map $l : X \to \C$), a
complete linear flag $\{ p \} = F_0 \subset F_1 \subset \dots
\subset F_d = X$, a closed ball $B \subset X$ around $p$, and an
algebraic Whitney stratification $\cX$ refining $\cS$, such that
the following three conditions hold.

(i)   For $S \in \cX$, if $p \notin \bar S$ then $\bar S \cap B =
\emptyset$.

(ii)  Let $Q_i = X / F_{d-i}$ and $\pi_i : V \to Q_i$ be the
projection.  Then for every $S \in \cX$ of dimension $i > 0$, the
map $\pi_{i - 1} \oplus l : S \cap B \to Q_{i - 1} \oplus \C$
has full rank.

(iii) Let $X_i$ be the union of all $S \in \cX$ with $\dim S
\leq i$.  Then for every $i = 1, \dots, d$, the intersection
$X_i \cap F_{d - i + 1} \cap \ker l \cap B = \{ p \}$.
\end{lemma}

Assuming Lemma \ref{flag} for the moment, we take $f$ in Theorem
\ref{main} to be the real part of the covector $l$, and we let the
ball $B$ to be the same as in the lemma.
We now describe the subset $K \subset B$.  It is defined as a union
$K = K^+ \cup K^-$, where $K^+$ and $K^-$ are constructed inductively.
Let $K^+_0 = K^-_0 = \{ p \}$.  Suppose now $1 \leq i \leq d$, and the
sets $K^\pm_{i-1}$ have been constructed.  We set
$$K^+_i = \{x \in X_i \cap B \,|\, \im l (x) = 0 \;\&\;\exists\; y \in
K^+_{i-1}: \pi_{i - 1} (x) = \pi_{i - 1} (y), \, f(x) \geq f(y)\},$$
$$K^-_i = \{x \in X_i \cap B \,|\, \im l (x) = 0 \;\&\;\exists\; y \in
K^-_{i-1}: \pi_{i - 1} (x) = \pi_{i - 1} (y), \, f(x) \leq f(y)\},$$
where $\im l : X \to \R$ is the imaginary part of $l$.  Lastly, we set
$K^\pm = K^\pm_d$.  It is clear that $K$ is a closed, real
semi-algebraic subset of $B$.

\begin{lemma}  With these definitions, conditions (i) and (ii) of
Theorem \ref{main} are satisfied.
\end{lemma}

\noindent
{\bf Proof:}  To check condition (i), note that $K^\pm \cap X_i =
K^\pm_i$.  We now prove by induction on $i$ that $\dim^{}_\R K^\pm_i
\leq i$.  Indeed, case $i = 0$ is trivial, and the induction step
follows from the definition of $K^\pm_i$ and condition (ii) of Lemma
\ref{flag}.

To check condition (ii), we argue by contradiction.  Suppose the set
$K^+ \cap f^{-1} (0)$ is larger than $\{ p \}$ (the case of $K^-$ is
of course analogous).  Let $i$ be the smallest integer such that
$K^+_i \cap f^{-1} (0)$ contains a point $x \neq p$.  By construction,
there is a $y \in K^+_{i-1}$ with $\pi_{i - 1} (x) = \pi_{i - 1} (y)$
and $f(x) \geq f(y)$.  By the minimality of $i$, we must have $y = p$
and, therefore, $x \in F_{n - i + 1}$.  But $x$ is also in $\ker l$,
since $\re l(x) = 0$ by assumption and $\im l (x) = 0$ by the
construction of $K^+_i$.  Thus we have a contradiction with condition
(iii) of Lemma \ref{flag}.
\hfill$\Box$

\vspace{.1in}

\noindent
{\bf Proof of Lemma \ref{flag}:}  Without loss of generality, we can
assume that every stratum of $\cS$ is connected (i.e., irreducible).
All stratifications in this proof will have connected strata, so we
can refer to ``the generic point of a stratum'' with no ambiguity.
We proceed inductively, starting with $i = 1$ and going up to $i = d$,
to construct the following four things:

$\;\, \bullet \;$  the $i$-plane $F_i \subset V$;

$\;\, \bullet \;$  the set $X_{d-i} \subset X$;

$\;\, \bullet \;$  a covector $l_i \in \Delta$;

$\;\, \bullet \;$  an (algebraic, Whitney) refinement $\cX_i$ of the
stratification $\cS$.

\noindent
When the process is complete, we will take $l = l_d$ and $\cX = \cX_d$.
After the $i$-th step of the construction, the following conditions will
be satisfied:

(1)  $X_{d-i}$ is the union of $S \in \cX_i$ with $\dim S \leq d - i$;

(2)  $\dim X_{d-i} \cap \ker l_i \leq d-i-1$;

(3)  if $i > 1$, every $S \in \cX_{i-1}$ with $\dim S \geq d - i + 2$ is
also a stratum of $\cX_i$;

(4)  if $i > 1$, we have $l_i|^{}_{F_i} = l_{i-1}|^{}_{F_i}$;

(5)  $p$ is an isolated point of $F_i \cap X_{d-i}$;

(6)  for every $S \in \cX_i$ with $\dim S = d - i$ and $p \in \bar S$,
the projection $\pi_{d-i} : S \to Q_{d-i}$ \\
has full rank near $p$;

(7)  $p$ is an isolated point of $F_i \cap X_{d-i+1} \cap \ker l_i$;

(8)  for every $S \in \cX_i$ with $\dim S = d - i + 1$ and
$p \in \bar S$, the map \\
$\pi_{d-i} \oplus l_i: S \to Q_i \oplus \C$ has full rank near $p$.

\noindent
Note that (3) and (4) ensure that (7) and (8) will continue to hold if
we replace $l_i$ by $l_d$ and $\cX_i$ by $\cX_d$.  Thus, our
construction will prove the lemma.

As a base step of the induction, we take $\cX_1 = \cS$, so $X_{d - 1}$
is the union of all but the largest stratum of $\cS$.  Let $C(X_{d-1})$
be the normal cone of $X_{d - 1}$ at $p$.  We take $F_1$ to be any line
not contained in $C(X_{d-1})$, and $l_1$ to be any covector in $\Delta$
which does not vanish on $F_1$.  Conditions (1)--(8) for $i = 1$ are
clearly satisfied.

Assume now $i > 1$, and the first $i-1$ steps of the construction have
been completed.  To select the plane $F_i$, we consider two cones in
$Q_{d-i+1} = X / F_{i-1}$.  First, let $X'_{d-i}$ be the union of all
$S \in \cX_{i-1}$ with $\dim S \leq d - i$, and let $C_1$ be the normal
cone of the image $\pi_{d-i+1} (X'_{d-i})$ at $p$.  It is a proper,
closed cone in $Q_{d-i+1}$.  Second, let $C_2$ be the normal cone at $p$
of $\pi_{d-i+1} (\ker l_{i-1} \cap X_{d-i+1})$.  By part (2) of the
induction hypothesis, this too is a proper, closed cone in $Q_{d-i+1}$.
We now choose any line $L \subset Q_{d-i+1}$ not contained in
$C_1 \cup C_2$, and set $F_i = \pi_{d-i+1} (L) \subset V$.  This also
defines the projection $\pi_{d-i} : V \to Q_{d-i}$.  Let $\Sigma$ be the
set of all $S \in \cX_{i-1}$ with $\dim S = d - i + 1$ and $p \in \bar S$.

\vspace{.1in}

\noindent
{\bf Claim 1:}  For every $S \in \Sigma$, the map
$\pi_{d-i} \oplus l_{i-1}: S \to Q_{d-i} \oplus \C$
has full rank at the generic point of $S$.

\vspace{.1in}

\noindent
{\bf Proof:}  Suppose the claim is false.  Consider the intersection
$\Gamma = F_i \cap S$.  By part (6) of the induction hypothesis,
there is an open neighborhood $U \subset X$ of $p$, such that $\Gamma
\cap U$ is a smooth curve, cut out transversely as the zero set of
$\pi_{d-i}|^{}_{S \cap U}$.  By part (5) of the induction hypothesis,
and because the line $L$ in the construction of $F_i$ was chosen not to
lie in the cone $C_1$, we have $p \in \bar \Gamma$.  Further, since $L$
was also chosen not to lie in the cone $C_2$, we can assume that there
are no critical points of $l_{i-1}|^{}_\Gamma$ in $U$.  This means that
the map $\pi_{d-i} \oplus l_{i-1}: S \to Q_{d-i} \oplus \C$ has full
rank at every point of $\Gamma \cap U$.
\hfill$\Box$

\vspace{.1in}

For $S \in \Sigma$, let $S^\circ$ be the part of $S$ where the map
$\pi_{d-i} \oplus l_{i-1}: S \to Q_{d-i} \oplus \C$ has full rank.
We set
$$X_{d-i} = X'_{d-i} \cup \bigcup_{S \in \Sigma} S \setminus S^\circ.$$
It is clear from the proof of Claim 1 that condition (5) is satisfied.

\vspace{.1in}

\noindent
{\bf Claim 2:}  Let $S \in \Sigma$, and let $T \subset S \setminus
S^\circ$ be an irreducible, smooth, locally closed subvariety with
$\dim T = d - i$ and $p \in \bar T$.  Then the restriction
$\pi_{d-i}|^{}_T$ has full rank at the generic point of $T$.

\vspace{.1in}

\noindent
{\bf Proof:}  Suppose the claim is false.  Then we have
$\dim \pi_{d-i} (\bar T) < d-i$.  Therefore, every fiber of
$\pi_{d-i}|^{}_{\bar T}$ must have positive dimension at every point.
But it is clear from the proof of Claim 1 that $p$ is an isolated
point of $F_i \cap \bar T$.
\hfill$\Box$

\vspace{.1in}

We are now ready to describe the stratification $\cX_i$.  We take every
stratum $S \in \cX_{i-1}$ with $\dim S > d-i$ and $S \notin \Sigma$ to
be also a stratum of $\cX_i$.  This ensures that conditions (1) and (3)
are satisfied.  For $S \in \Sigma$, we take each irreducible component
of $S^\circ$ to be a stratum of $\cX_i$.  It remains to stratify the
set $X_{d-i}$.  By Claim 2, this can be done in such a way that the
resulting stratification $\cX_i$ is a Whitney refinement of $\cX_{i-1}$,
and condition (6) is satisfied.

The last thing to construct is the covector $l_i$.  To satisfy condition
(2) we must ensure that $l_i$ does not vanish identically on any of the
$(d-i)$-dimensional strata of $\cX_i$.  We take $l_i = l_{i-1} + h \circ
\pi_{d-i}$, where $h : Q_{d-i} \to \C$ is a small linear functional in
general position.  It is easy to check using condition (5) that, for a
suitable choice of $h$, we will have $l_i \in \Delta$ and condition (2)
will hold.

This completes the construction of the quadruple $\{ F_i, X_{d-i}, l_i,
\cX_i\}$.  We have already remarked that conditions (1), (2), (3), (5),
and (6) are satisfied.  Condition (4) is clear from the definition of
$l_i$.  Condition (7) follows from the fact that the line $L$ in the
construction of $F_i$ was chosen not to lie in the cone $C_2$, combined
with condition (4) and part (5) of the induction hypothesis.  Finally,
condition (8) follows from the definition of the loci $S^\circ$
($S \in \Sigma$), again combined with condition (4).
\hfill$\Box$

\section{Construction of the Flows}

In this section, we complete the proof of Theorem \ref{main} by
constructing the $\nabla f$-like vector fields appearing in part
(iii) of that theorem.  We begin with some preparations (keeping
the notation of Section 4).

Theorem \ref{main} stipulates that the vector field $V$ should be
$\cS$-preserving, i.e., weakly controlled with respect to some
control data on $(X, \cS)$.  However, by Proposition \ref{approx},
it suffices to construct, for each $\cU \supset K$, a vector field
$V$ which is $\cX$-preserving, instead.  We therefore fix some
control data on $(X, \cX)$, subject to the only condition that the
quasi-distance function $\rho_{\{ p \}}$ is the standard Euclidean
distance to $p$.  All weakly controlled vector fields in the rest
of this section will be compatible with this control data.

Let $\cX'$ be the set of all $S \in \cX$ with $S \cap B \neq
\emptyset$ and $\dim S > 0$.  For every $S \in \cX'$, with
$\dim S = i$, let $V^S$ be the unique smooth vector field on
$S \cap B$ satisfying the following properties:

(a)  the (standard, Euclidean) norm $|\!| V^S_x |\!| = 
|\!| x |\!|$ for every $x \in S \cap B$;

(b)  $V^S$ preserves the projection $\pi_{i-1} :
S \cap B \to Q_{i-1}$;

(c)  $V^S$ preserves the imaginary part $\im l :
S \cap B \to \R$;

(d)  $V^S_x f > 0$ for every $x \in S \cap B$.

\noindent
For every $S \in \cX'$, we fix a continuous, controlled extension
$\tilde V^S$ of $V^S$ to some neighborhood $U_S$ of $S \cap B$.
The vector field $V$ 
will be constructed by ``patching together'' the $\tilde V^S$.

Without loss of generality, we may assume that $U_S \cap T =
\emptyset$, for all distinct $S, T \in \cX'$ with
$\dim T \leq \dim S$.  Let $\tilde U_S$ be the union of all
$T \in \cX'$ with $\dim T > \dim S$ or $T = S$.  We say that
$\phi : \tilde U_S \to [0,1]$ is a cut-off function for
$S \cap B$ if $\supp (\phi) \subset U_S$ and $\phi^{-1} (1)$
contains a neighborhood of $S \cap B$.  Given two cut-off
functions $\phi, \psi$ for $S \cap B$, we write $\phi \less \psi$
if $\supp(\psi) \subset \phi^{-1}(1)$.

Now, let $\Phi$ be the set of all collections $\bar \phi =
\{ \phi_S  \}_{S \in \cX'}$, where $\phi_S$ is a cut-off function
for $S \cap B$.  Given $\bar \phi, \bar \psi \in \Phi$, we write
$\bar\phi \less \bar\psi$ if $\phi_S \less \psi_S$ for all $S \in
\cX'$.  It is easy to see that $(\Phi, \less)$ is a directed set:
$$\forall \; \bar\phi_1, \bar\phi_2 \in \Phi \;\; \exists \;
\bar\phi_3 \in \Phi : \; \bar\phi_1 \less \bar\phi_3, \; 
\bar\phi_2 \less \bar\phi_3.$$

Fix an ordering $\cX' = \{ S_1, \dots, S_n \}$ such that
$\dim S_i \leq \dim S_j$ for $i \leq j$.  For every $\bar \phi \in
\Phi$, we define a vector field $V = V(\bar\phi)$ on $B$ as follows.
Let $V_p = 0$, and use the formula
$$V = \phi_{S_1} \cdot \tilde V^{S_1} + (1 - \phi_{S_1}) \cdot (
\phi_{S_2} \cdot \tilde V^{S_2} + (1 - \phi_{S_2}) \cdot (
\phi_{S_3} \cdot \tilde V^{S_3} + \dots +
(1 - \phi_{S_{n-1}}) \cdot V^{S_n}) \dots )$$
on $B \setminus \{ p \}$.
This formula should be parsed left to right, and evaluation should
stop as soon as one of the expressions $(1 - \phi_{S_i})$ is found
to be zero.  In this way, we will never have to evaluate one of
the $\phi_{S_j}$ or $\tilde V^{S_k}$ at a point where it is not
defined.  Also, note the missing $\;\tilde{}\;$ over $V^{S_n}$; it
is not needed because $S_n$ is the open stratum of $\cX$.  It is
not hard to check that there is a $\bar \phi_0 \in \Phi$ such that,
for every $\bar \phi \in \Phi$ with $\bar\phi_0 \less \bar\phi$, the
vector field $V(\phi)$ is $\nabla f$-like.  Theorem \ref{main} is a
consequence of the following claim.

\vspace{.1in}

\noindent
{\bf Claim 1:}  For every open $\cU \supset K$, there exists a
$\bar \phi_0 \in \Phi$ such that, for every $\bar \phi \in \Phi$
with $\bar\phi_0 \less \bar\phi$, we have
$M^\pm_{V(\phi)}(p) \subset \cU$.

\vspace{.1in}

Claim 1, in turn, follows from the following.

\vspace{.1in}

\noindent
{\bf Claim 2:}  For every $x \in B \setminus K$, there is a
neighborhood $U_x$ of $x$ and a $\bar \phi_0 \in \Phi$ such that,
for every $\bar \phi \in \Phi$ with $\bar\phi_0 \less \bar\phi$,
we have $M^\pm_{V(\phi)}(p) \cap U_x = \emptyset$.

\vspace{.1in}

Claim 2 is readily proved by induction on the dimension of the
stratum containing $x$, using the definition of $K$ and the
following.

\vspace{.1in}

\noindent
{\bf Claim 3:}  Fix an $i \in \{ 0, \dots, d-1 \}$, and let
$A \subset B$ be a compact set such that $A \cap X_i = \emptyset$.
Then for every $\epsilon > 0$, there is a $\bar \phi_0 \in \Phi$
such that, for every $\bar \phi \in \Phi$ with $\bar\phi_0 \less
\bar\phi$ and every $a \in A$, we have the following estimates:
$$|\!| d \pi_i (V(\phi)_a) |\!| < \epsilon \cdot d f (V(\phi)_a),$$
$$| \im d l (V(\phi)_a) | < \epsilon \cdot d f (V(\phi)_a).$$

\vspace{.1in}

Lastly, Claim 3 follows from the fact that the $\tilde V^S$ in the
construction of $V(\phi)$ were chosen to be continuous extensions
of the vector fields $V^S$ satisfying conditions (b)-(d).  This
completes the proof of Theorem \ref{main}.

\vfill

\section{Self-Indexing}

In this section, we give a proof of Theorem \ref{siinsm}.  It
follows the same scheme as the proof of Theorem \ref{siincm}
outlined in Section 1.1.

\begin{defn}\label{sss}
Let $(X, \cS)$ be a compact, Whitney stratified $C^\infty$
manifold.  A stratified subset $A \subset X$ is a closed
subset presented as a finite disjoint union $A = \bigcup A_i$
so that each $A_i$ is a smooth submanifold of one of the strata
of $\cS$, and the frontier $\overline{A_i} \setminus A_i$ is a
union of several of the $A_j$ with $\dim A_j < \dim A_i$.
\end{defn}

\begin{defn}\label{tdep}
Let $(X, \cS)$ be a compact, Whitney stratified $C^\infty$
manifold, with fixed control data.  A time-dependent
controlled vector field $\{ V_t \}_{t \in [0,1]}$ on $(X, \cS)$
is a controlled vector field defined in some neighborhood of
$X \times [0,1] \subset X \times \R$, whose component in
the $\R$-direction is identically zero.
\end{defn}

\begin{lemma}\label{genpos}
Let $(X, \cS)$ be as in Definition \ref{tdep}, and let $A, B
\subset X$ be two stratified subsets.  Assume that for every
$S \in \cS$, we have:
$$\dim (A \cap S) + \dim (B \cap S) < \dim S.$$
Then there exists a time-dependent
controlled vector field $\{ V_t \}_{t \in [0,1]}$ on $X$ whose
time-1 flow $\psi_{V,1} : X \to X$ satisfies $\psi_{V,1} (A)
\cap B = \emptyset$.
\end{lemma}

\noindent
{\bf Proof:}  The restriction of $V_t$ to the $i$-skeleton of
$\cS$ is constructed by induction on $i$.  The induction step
is an application of the general position in manifolds.  See
[12, \S3] for a proof of a much more general stratified
general position result.
\hfill$\Box$

\vspace{.1in}

Theorem \ref{siinsm} follows from Lemmas \ref{interchange} 
and \ref{allign} below.

\begin{lemma}\label{interchange}
Let $(X,\cS)$ be a Whitney stratified complex algebraic
variety, and let $g : X \to \R$ be a proper Morse function.
Let $I \subset \R$ be an open interval whose preimage contains
exactly two critical points: $g^{-1}(I) \cap \Sigma_f = \{p,q\}$.
Assume that $\index_g(p) \leq \index_g(q)$ and $g(p) \geq g(q)$.
Then for every $a,b \in I$, there is a Morse function
$f : X \to \R$ such that:

(i)   $f^{-1}(I) = g^{-1}(I)$, and $f = g$ outside of some
compact subset of $g^{-1}(I)$;

(ii)  $f^{-1}(I) \cap \Sigma_f = \{p,q\}$;

(iii) in some neighborhood of $p$ we have 
$f(x) = g(x) + a - g(p)$;

(iv)  in some neighborhood of $q$ we have 
$f(x) = g(x) + b - g(q)$.
\end{lemma}

\begin{lemma}\label{allign}
Let $(X, \cS)$ and $g : X \to \R$ be as in Lemma
\ref{interchange}, and let $I \subset \R$ be an open interval
such that all the critical points in $g^{-1}(I)$ have the same
index.  Then for every $c \in I$, there is a Morse function
$f : X \to \R$ such that:

(i)   $f^{-1}(I) = g^{-1}(I)$, and $f = g$ outside of some
compact subset of $g^{-1}(I)$;

(ii)  $f^{-1}(I) \cap \Sigma_f = g^{-1}(I) \cap \Sigma_g$;

(iii) near each $p \in f^{-1}(I) \cap \Sigma_f$ we have
$f(x) = g(x) + c - g(p)$.
\end{lemma}

The proofs of Lemmas \ref{interchange} and \ref{allign} are
similar; we will only give the first.

\vspace{.1in}

\noindent
{\bf Proof of Lemma \ref{interchange}:}
{\it Step 1:}  The case when $g(p) = g(q)$ is obvious, so we
assume that $g(p) > g(q)$.  Pick two numbers $c < d$ from the
interval $(g(q), g(p))$.  Let $D = g^{-1}[c,d]$.  Fix control
data $\{ U_S, \Pi_S, \rho_S \}$ on $(X, \cS)$ which is
$g$-compatible in some neighborhood of $D$.

\vspace{.1in}

\noindent
{\bf Claim 1:}  There exist a $\nabla g$-like vector field
$V^D$ on $D$, compatible with $\{ U_S, \Pi_S, \rho_S \}$, and a
pair of stratified subsets $A, B \subset g^{-1}(d)$ such that:

(i)  for every $S \in \cS$ with $\dim_\C S = s$, we have:
$$\dim_\R (A \cap S) < s - \index_g(q);$$
$$\dim_\R (B \cap S) < s + \index_g(p).$$

(ii) for every pair of open neighborhoods $\cU_A, \cU_B \subset
g^{-1}(d)$ of $A$ and $B$, there is a $\nabla g$-like vector
field $V$ on $X$ with $V |^{}_D = V^D$ such that:
$$(M^+_V(q) \cap g^{-1}(d)) \subset \cU_A;$$
$$(M^-_V(p) \cap g^{-1}(d)) \subset \cU_B.$$

\vspace{.1in}

Claim 1 follows by putting together Theorem \ref{main},
Proposition \ref{approx}, and Corollaries \ref{ml} and \ref{ns}.

{\it Step 2:}  Using Lemma \ref{genpos} to modify the vector
field $V^D$, we can strengthen Claim 1 as follows.

\vspace{.1in}

\noindent
{\bf Claim 2:}  The sets $A$ and $B$ in Claim 1 can be chosen
to satisfy $A \cap B = \emptyset$. 

\vspace{.1in}

This immediately implies that there exists a $\nabla g$-like
vector field $V$ on $X$, compatible with $\{ U_S, \Pi_S, 
\rho_S \}$, such that $M^-_V (p) \cap M^+_V (q) = \emptyset$.

{\it Step 3:}  Now, choose a smooth function $\phi : g^{-1}(d)
\to [0,1]$ such that:

(a)  the control data $\{ U_S, \Pi_S, \rho_S \}$ on $(X, \cS)$
restricts to $\phi$-compatible control data on $g^{-1} (d)$
(with the stratification induced from $\cS$);

(b)  there is a neighborhood $\cU_A \subset g^{-1}(d)$ of
$M^+_V(q) \cap g^{-1}(d)$ such that $\phi|^{}_{\cU_A} = 0$;

(c)  there is a neighborhood $\cU_B \subset g^{-1}(d)$ of
$M^-_V(p) \cap g^{-1}(d)$ such that $\phi|^{}_{\cU_b} = 1$.

\noindent
It is not hard to check that $\phi$ extends uniquely to a smooth
function $\tilde \phi : g^{-1}(I) \to [0,1]$, which is constant
along the flow lines of $V$.

The function $f$ is now easy to construct by setting
$f|^{}_{g^{-1}(I)} = \tilde f \circ \alpha$, where
$\alpha = (g, \phi) : g^{-1}(I) \to I \times [0,1]$
and $\tilde f : I \times [0,1] \to I$ is a suitable smooth
function of two variables.
\hfill$\Box$

\vspace{.1in}

This completes the proof of Theorem \ref{siinsm}.

\vspace{.1in}

\noindent
Northwestern University, Department of Mathematics,
2033 Sheridan Road \\
Evanston, IL 60208-2730 \hfill
{\em grinberg@math.nwu.edu}

\end{document}